\documentclass[11pt]{amsart}
\usepackage{times}      % use Times fonts for text
\usepackage{latexsym}   % use all LaTeX fonts
\usepackage{amssymb}    % use all AMS fonts
\usepackage{amsmath}    % include some AMS-LaTeX functionality
\usepackage{amsthm}     % AMS style theorem and proof environments
\usepackage{mathrsfs}   % special calligraphic characters like \mathscr{S}
\usepackage{euscript}   % euscript package
\usepackage{epsfig}             % epsfig package, encap postscript
\usepackage[all]{xy}    % Xy-Pic for diagrams
\usepackage{amsmath,amsthm,amssymb}
\input xypic
\usepackage{graphicx}

\theoremstyle{plain}

\theoremstyle{remark}

\theoremstyle{definition}

\def\Box{\square}
\def\ptorus{\overset{\circ}{T}}

\def\ptorusn{\ptorus\times\cdots\times\ptorus}

\def\textdiam{\,\hbox{diam\,}}

\def\textindsmall{\hbox{\,Ind\,}}
\def\textind{\hbox{\,Ind}}

\def\textdim{\mathrm{dim\,}}

\def\calv{{\mathcal V}}
\def\cala{{\mathcal A}}

\def\calm{{\mathcal M}}
\def\calb{{\mathcal B}}
\def\calc{{\mathcal C}}

\def\calh{{\mathcal H}}

\def\cgm{C^\ast_G(M)}
\def\dgm{D^\ast_G(M)}

\def\reals{{\mathbb{R}}}
\def\rpn{\reals{\mathbb{P}}^n}

\def\integers{{\mathbb{Z}}}
\def\textsupp{\hbox{Supp}\,}

\def\textdiag{\hbox{\,\,diag\,}}

\def\textaut{\hbox{\,Aut\,}}

\def\textpen{\hbox{Pen}}

\def\slashs{\,\hbox{\slash}\kern-6.5pt S}
\def\slashd{\,\hbox{\slash}\kern-8.0pt D}

\def\rationals{{\mathbb Q}}
\def\calr{{\mathcal R}}

\def\cald{{\mathcal D}}
\def\calds{\,{{\slash}\kern-8.0pt \cald }}

\def\frakp{{\mathfrak p}}
\def\fraka{{\mathfrak a}}
\def\frakg{{\mathfrak g}}
\def\frakf{{\mathfrak f}}
\def\slnr{\mathrm{SL}_n(\reals)}
\def\slnz{\mathrm{SL}_n(\integers)}
\def\sonr{\mathrm{SO}_n(\reals)}

\begin{document}

\title{Coarse Obstructions to Positive Scalar Curvature Metrics
in Noncompact Arithmetic Manifolds}
\author{Stanley S. Chang}
\address{Rice University, Houston, TX 77005}
\email{sschang@math.rice.edu}
% \keywords{KEYWORDS}

\begin{abstract}
Block and Weinberger show that an arithmetic manifold
can be endowed with a positive scalar curvature metric if and only
if its $\rationals$-rank exceeds $2$. We show in this article that
these metrics are never in the same coarse class as the natural
metric inherited from the base Lie group. Furthering the coarse
$C^\ast$-algebraic methods of Roe, we find a nonzero Dirac obstruction
in the $K$-theory of a particular operator algebra which encodes
information about the quasi-isometry type of the manifold as well
as its local geometry.

\end{abstract}

\maketitle

\bigskip\noindent
I. Introduction

\bigskip
In the course of  showing
 that  no manifold of non-positive sectional curvature
can be endowed with a metric of positive scalar curvature,
Gromov and Lawson \cite{GromovLawson2} were led to  consider
what we would now call
restrictions on the coarse equivalence type of complete noncompact
manifolds of such positively curved metrics. In particular, they showed
that such metrics cannot exist in manifolds 
for which there exists a degree one proper Lipschitz map from
the universal cover to $\reals^n$, now understood to be
essentially a coarse condition. 
Block and Weinberger \cite{BlockWeinberger} investigate
the situation in which no coarse conditions are imposed upon
the complete metric, 
focusing on quotients
$\Gamma\backslash G/K$ of symmetric spaces associated to a lattice
$\Gamma$ in an irreducible semisimple Lie group $G$. They show that
the space $M=\Gamma\backslash G/K$ can be
given a complete metric of uniformly positive scalar curvature 
$\kappa\ge\varepsilon >0$ if and only if $\Gamma$ is an arithmetic
group of $\rationals$-rank exceeding $2$. 

Note that the theorem of Gromov and Lawson \cite{GromovLawson2}
mentioned above establishes this theorem in the case of
$\hbox{rank}_\rationals\Gamma=0$. 
In the higher rank cases, for which the resulting quotient
space is noncompact, the metrics constructed by Block and Weinberger
are however wildly different in the large when compared to the natural one
on $M$ inherited from the base Lie group $G$. In fact, their examples  are
all coarse quasi-isometric to rays. Their theory evokes
a natural question: Can the metric be chosen so that it is
simultaneously uniformly positively curved and coarsely equivalent to the
natural metric induced by $G$? 

One of the important developments in analyzing positive scalar
curvature in the context of noncompact manifolds, especially 
when restricted to the coarse quasi-isometry type, 
is introduced by 
Roe \cite{RoeIndex}, \cite{RoeCoarse}, who considers a higher
index, analogous to the Novikov higher signature, that lives naturally
in the $K$-theory of the $C^\ast$-algebra $C^\ast(M)$ of operators on $M$ with
finite propagation speed. 
He describes a map from the $K$-theory group $K_\ast(C^\ast(M))$
to the $K$-homology $K_\ast(\nu M)$ of the Higson corona space
which admits a dual transgression map 
$H^\ast(\nu M)\to HX^\ast(M)$.  If the Dirac operator on $M$ is invertible, then
the image of its index in $K_\ast(\nu M)$ vanishes, leading to
vanishing theorems for the index paired with coarse classes from
the transgression of $\nu M$.  Roe's construction is used to show
that a metric on a noncompact manifold
cannot be uniformly positively curved if the Higson corona of the manifold
contains an essential $(n-1)$-sphere.  Such spaces are called
{\it ultraspherical manifolds}.

The usual Roe algebra, however, is unsuited to provide information
about the existence of positive scalar curvature metrics
that exist on arithmetic manifolds, because their coronae
are too anemic. For example,
the space at infinity of a product of punctured two-dimensional 
tori is a simplex and therefore contractible.
As a coarse object, the $K$-theory of the
Roe algebra associated to this multi-product 
space can be identified
with $K_\ast(C^\ast(\reals_{\ge0}^n))$.
Yet Higson, Roe and Yu \cite{HigsonRoeYu} have shown that the Euclidean
cone $cP$ on a single simplex $P$ must satisfy 
$K_\ast(C^\ast(cP))=0$. Since the Euclidean hyperoctant
$\reals^n_{\ge0}$ is simply the cone on an $(n-1)$-simplex,
we find that $K_\ast(C^\ast(\reals_{\ge0}^n))$
is the trivial group and hence no obstructions are detectable.
Even by considering the fundamental group of the manifold by tensoring
the Roe algebra with $C^\ast\pi_1(M)$ is this detection process unfruitful,
since the $K$-theory group $K_\ast(C^\ast(M)\otimes C^\ast\pi_1(M))$
vanishes as well. What seems to be critical is how different elements
of the fundamental group at infinity can be localized to different parts
of the space at infinity.

In this article, we shall
 provide coarse indicial obstructions in the following noncompact
manifolds: a finite product of punctured two-dimensional tori, a
finite product of hyperbolic manifolds, the double quotient
space $\slnz\backslash \slnr /\sonr$ of unit volume tori, and
more generally the double quotient space $\Gamma\backslash G/K$,
where $G$ is an irreducible semisimple Lie group, $K$ its maximal
compact subgroup and $\Gamma$ an arithmetic subgroup of $G$.
Note that the first two do not correspond to irreducible
quotients, but an analysis of these spaces gives us the proper
insight to attack the more general cases. A further research project
will analyze this problem without the irreducibility assumption.
The key feature in these particular manifolds $M$ is that they contain
hypersurfaces $V$ that are coarsely equivalent to a product $E\times U$
of Euclidean space $E$ with some iterated circle bundle $U$ (i.e. a torus,
Heisenberg group, or more generally a group of unipotent matrices).
Moreover such a hypersurface decomposes the manifold $M$ into a
{\it coarsely excisive pair} $(A,B)$ for which $A\cup B=M$ and
$A\cap B=V$. A generalized form of the Mayer-Vietoris sequence
constructed by Higson, Roe and Yu \cite{HigsonRoeYu}
provides the following:
$$\cdots\longrightarrow K_\ast(C^\ast_G(A))\oplus K_\ast(C^\ast_G(B))
\longrightarrow K_\ast(C^\ast_G(M))\longrightarrow K_{\ast-1}(C^\ast_G(V))
\longrightarrow\cdots$$
The boundary map $\partial: K_\ast(C^\ast_G(M))
\longrightarrow K_{\ast-1}(C^\ast_G(E\times U))$ sends
$\textindsmall_M(D)$, the
index of the spinor Dirac bundle on the universal cover lifted from
that on $M$, to $\textindsmall_{E\times U}(D)$. To see
that these indices are indeed nonzero, we note that there is a
boundary map $K_{\ast-1}(C^\ast_G(E\times U))\to K_{\ast-\textdim E}
(C^\ast_G(\reals\times U))$, 
which sends index to index.  We show that the index of the Dirac
operator in the latter
group,  however, is nonzero by noting that the Gromov-Lawson-Rosenberg
conjecture is true for nilpotent groups and hence provides an appropriate
nonzero obstruction.

I would like to thank
Alex Eskin, Benson Farb, Nigel Higson, Thomas Nevins,
Mel Rothenberg, John Roe, Stephan Stolz and Guoliang Yu for very
useful conversations. In particular,
I would like to acknowlege the role of 
my advisor Shmuel Weinberger in pointing out the strength of certain
tools in the realization of these theorems.

\bigskip\bigskip\noindent
II. The Generalized Roe Algebra

\bigskip
The {\it coarse category\,} is defined to contain metric spaces as
its objects and maps $f: (X,d_X)\to (Y,d_Y)$
between metric spaces as its morphisms satisfying
the following expansion and properness conditions: (a) for each
$R>0$ there is a corresponding $S>0$ such that, if $d_X(x_1,x_2)\le R$
in $X$, then $d_Y(f(x_1),f(x_2))\le S$, (b) the inverse image 
$f^{-1}(B)$ under $f$
of each bounded set $B\subseteq Y$ is also bounded in $X$. Such a
function will be designated a {\it coarse map}, and two coarse maps
$f,g: X\to Y$ are said to be {\it coarsely equivalent} if their mutual
distance of separation
$d_Y(f(x), g(x))$ is uniformly bounded in $x$. Naturally
two metric spaces are coarsely equivalent if there exist maps from
one to the other whose compositions are coarsely equivalent to the
appropriate identity maps. Two metrics $g_1$ and $g_2$ 
on the same space $M$ are said to be coarsely equivalent if 
$(M,g_1)$ and $(M,g_2)$ are coarsely equivalent metric spaces.

\bigskip
Following Roe \cite{RoeCoarse}, we recall a Hilbert space $H$
is an $M$-module for a manifold $M$ if there is a representation
of $C_0(M)$ on $H$, that is, a $C^\ast$-homomorphism
$C_0(M)\to B(H)$.  We will say that an operator $T : H\to H$ is
{\it locally compact} if, for all $\varphi\in C_0(M)$, the operators
$T\varphi$ and $\varphi T$ are compact on $H$. We define the
{\it support} of $\varphi$ in an $M$-module $H$
to be the smallest closed set $K\subseteq M$
such that, if $f\in C_0(M)$ and $f\varphi\ne0$, then $f\vert_K$
is not identically zero. Consider the 
$\widetilde M$-module $H=L^2(\widetilde M)$, where $\widetilde M$
is the universal cover of $M$ endowed with the appropriate metric lifted
from the base space. Let $\pi: \widetilde M\to M$ be the usual projection
map and for any $\varphi, \psi\in C_0(\widetilde M)$ consider the collection 
$\Gamma(\varphi, \psi)$ of paths $\gamma: [0,1]\to \widetilde M$
in $\widetilde M$ 
originating in $\textsupp (\varphi)$ and ending
in $\textsupp (\psi)$. Denote by $L[\,\gamma\,]$, 
for $\gamma\in\Gamma(\varphi, \psi)$,
the maximum distance of any two points on the projection of 
the curve $\gamma$ in $M$ by $\pi$, i.e.
$L[\,\gamma\,]= \sup_{x,y\in [0,1]}d(\pi\circ\gamma(x),
\pi\circ\gamma(y))$. 

\medskip\noindent
{\bf Definition:} Let $M$ be a manifold with universal cover 
$\widetilde M$. We say that an operator $T$ on  $L^2(\widetilde M)$
has {\it generalized finite propagation} if there is a constant $R>0$
such that $\varphi T\psi$ is identically zero in $B(H)$ whenever 
$\varphi,\psi\in C_0(\widetilde M)$ satisfies $$\inf_{\gamma\in
\Gamma(\varphi,\psi)}L[\,\gamma\,] >R.$$
The infimum of all such $R$ will be the {\it generalized propagation
speed} of the operator $T$. If $G=\pi_1(M)$ is the
fundamental group of $M$, we denote by 
$D^\ast_G(M)$ to be the norm closure of the $C^\ast$-algebra 
 of all locally compact, $G$-equivariant, generalized finite
propagation operators on $H$. 

\medskip
Let $M$ be a manifold and $\widetilde M$ its universal cover. Let
$T:H\to H$ be an operator on $H=L^2(\widetilde M)$.
Consider the subset $Q\subseteq \widetilde M\times\widetilde M$ of
pairs $(m,m')$ for which there exist functions $\varphi, \psi
\in C_0(\widetilde M)$ such that $\varphi(m)\ne0$, $\psi(m')\ne0$
and $\varphi T\psi$ does not identically vanish. We will say that
the {\it support} of $T$ is the complement in $\widetilde M\times
\widetilde M$ of $Q$. For such two points $m,m'\in\widetilde M$,
let $\gamma_{mm'}: [\,0,1\,]\to \widetilde M$
be the path of least length joining $m$ and $m'$ in $\widetilde M$.
We consider the projection of this path into $M$ by $\pi$ and
take the greatest distance between two points on this projected path.
Then it is easy to see that an operator $T$
has generalized finite propagation, as previously defined,  if
$$\sup_{m,m'}\,\,
\sup_{x,y\in[0,1]}d(\pi\circ\gamma_{mm'}(x), \pi\circ\gamma_{mm'}(y))
<\infty.$$

\medskip\noindent
{\bf Definition:}  Consider the norm closure $I$ of the
ideal in $D^\ast_G(M)$ generated by operators $T$
whose matrix representation,  parametrized by $\widetilde M
\times \widetilde M$, satisfies the
condition that  $(\pi\times\pi)(\textsupp T)$ is bounded
in $M\times M$.
Then the {\it generalized Roe algebra}, denoted
by $C^\ast_G(M)$, is obtained as
the quotient $D^\ast_G(M)/I$. Two operators in
$D^\ast_G(M)$ belong to the same class in $C^\ast_G(M)$
if their nonzero entries differ on at most
a bounded set when viewed from the perspective of
the base space.

\smallskip\noindent
{\bf Examples:} 

(1)  Let $T: L^2(\reals)\to L^2(\reals)$ be operator
on $L^2$-functions on the real line given by $(Tg)(x)=g(x+1)$ for
all $g\in L^2(\reals)$ and $x\in\reals$. Then for any $\varphi,\psi
\in C_0(\reals)$, $(\varphi T\psi)g (x) = \varphi(x+1)g(x+1)\psi(x)$.
If $\varphi$ is supported at $m=1$ and $\psi$ is supported at $m'=0$,
then $(\varphi T\psi)g$ is nonzero for any $g$ supported at $x=1$.
Hence $(0,1)\in\textsupp T$. It is easy to see that $(m,m')\in \textsupp T$
if and only if $m'-m=1$. The propagation speed of $T$ is $1$. If we write
$T$  as a matrix parametrized by $\reals\times\reals$,
all the nonzero entries
will lie at distance one from the diagonal.

(2) Let $M$ be the cylinder $S^1\times\reals$ with
its universal cover $\widetilde M=\reals^2$. An operator in the
algebra $D^\ast_G(M)$ will be some $T: H\to H$
on $L^2(\reals^2)$, which is 
of finite propagation speed (in the usual sense) in the direction
projecting down to the noncompact direction in $M$, but has no such
condition in the orthogonal direction corresponding to the compact
direction of $M$. In this direction, however, the operator is controlled
by the condition that it be $\integers$-equivariant. It is apparent
that the operator, when restricted to individual fibers, has finite
propagation speed, there is no requirement that 
the speed to be uniformly bounded across all fibers.

(3) Let $M=\,\,\stackrel{\circ}{\rpn}$, $n\ge3$, 
the once-punctured real projective
space, expressible as the quotient $(S^{n-1}\times\reals)/\integers_2$.
Certainly $M$ is coarsely equivalent to the ray $[\,0,\infty)$ and
is covered by the space $\widetilde M = S^{n-1}\times \reals$, where the points
$(s,r)$ and $(-s,-r)$ are identified by the projection map to $M$. Let
$T:L^2(\widetilde M)\to L^2(\widetilde M)$ be given by the reflection
$(Tf)(s,r)=f(s,-r)$.
Consider $\varphi_i,\psi_i\in C_0(\widetilde M)$ compactly supported
on $S^{n-1}\times [\,-i-1,-i\,]$ and $S^{n-1}\times [\,i,i+1\,]$,
respectively. Notice that $\varphi T\psi$ will never be identically zero,
and yet the length $L_i[\,\gamma\,]$ associated to $\varphi_i$ and
$\psi_i$ will always be at least $i$. Hence the operator $T$ 
is {\it not} of generalized  finite propagation speed 
and therefore not an element of the generalized Roe algebra $C_G^\ast(M)$.

\medskip\noindent
{\bf Lemma 1:}
Let $D$ a generalized elliptic operator in $L^2(M,S)$.
Suppose that $\widetilde D$ is the lifted operator on $\widetilde M$.
If $\Phi: \reals\to\reals$  is compactly supported,
then $\Phi(\widetilde D)$ lies in the generalized Roe
algebra $C^\ast_G(M)$.

\smallskip\noindent
{\it Proof:} (cf. \cite{RoeCoarse}, \cite{CheegerGromovTaylor})
Suppose that $\Phi$ has compactly supported Fourier
transform and denote by $\widehat\Phi$ the Fourier transform of $\Phi$.
 We may write $$\Phi(\widetilde D) =\frac{1}{2\pi} 
\int_{-\infty}^\infty\widehat\Phi(t)\,e^{it\widetilde D}\,dt.$$
It is known that $e^{it\widetilde D}$ has finite propagation speed,
and since $\widehat\Phi$ is compactly supported, the integral 
is defined and has a generalized propagation bound. Moreover,
by construction $\widetilde D$ is $\pi_1(M)$-equivariant. So
$\Phi(\widetilde D)$ is $\pi_1(M)$-equivariant as well. Therefore if
$\widehat \Phi$ is compactly supported, then 
$\Phi(\widetilde D)$ lies in $D^\ast_G(M)$ and passes to
an element of the quotient $C^\ast_G(M)$. However, functions
with compactly supported Fourier transform form a dense set in
$C_0(\reals)$ and the functional calculus map $f\mapsto f(\widetilde D)$
is continuous, so the result holds for all $\Phi\in C_0(\reals)$.
\hfill $\Box$

\medskip
Let $\chi:\reals\to\reals$ be a {\it chopping function} on $\reals$, i.e.
an odd continuous function with the property that $\chi(x)\to\pm 1$ as
$x\to\pm\infty$. In addition, denote by $B^\ast_G(M)$ the
{\it multiplier algebra} of $C^\ast_G(M)$, that is, the collection
of all operators $S$ such that $ST$ and $TS$ belong to $C^\ast_G(M)$
for all $T\in C^\ast_G(M)$. Then $B^\ast_G(M)$ contains $C^\ast_G(M)$
as an ideal. If $D$ is a generalized elliptic operator on $M$ and
$\widetilde D$ its lift to $\widetilde M$, then $\chi(\widetilde D)$ belongs
to $B^\ast_G(M)$. In addition, since $\chi^2-1\in C_0(\reals)$, 
we have $\chi(\widetilde D)^2-1\in C^\ast_G(M)$. Moreover, since
the $\integers_2$-grading renders the decompositions 
$$\chi(\widetilde D)=\left(\begin{array}{cc}
       0&\chi(\widetilde D)_-\\
       \chi(\widetilde D)_+&0
	\end{array}
  \right), ~~~~~
  \varepsilon=\left(\begin{array}{rr}
    1&0\\
    0&-1
    \end{array}
    \right), $$
it follows that $\varepsilon\,\chi(\widetilde D)+
\chi(\widetilde D)\,\varepsilon
=0$. By the discussion in \cite{RoeCoarse}, it follows
that  $F=\chi(\widetilde D)$ is a
Fredholm operator and admits an index $\textindsmall F
\in K_0(C^\ast_G(M))$.
In addition, any two chopping functions $\chi_1$
and $\chi_2$ differ by an element of $C_0(\reals)$. By the lemma
above, we have $\chi_1(\widetilde D)-\chi_2(\widetilde D)\in C^\ast_G(M)$,
so they define the same elements of $K$-theory. The common value 
for $\textindsmall F$ is denoted $\textindsmall
(D)$ and called the {\it generalized coarse index} of $D$. 
We write $C^\ast_G(M)$ and $\textindsmall (D)$ 
instead of $C^\ast_G(\widetilde M)$ and $\textindsmall(\widetilde D)$
to indicate that the construction is initiated by a generalized Dirac
operator on the base space. The following statements are standard
results of index theory; one may consult \cite{RoeCoarse} and
\cite{RoeIndex} for the essentially identical proof in the 
nonequivariant case.

\medskip\noindent
{\bf Proposition 1:}
Let $D$ be a generalized
elliptic operator in $L^2(M,S)$. 
If $0$ does not belong to
the spectrum of $\widetilde D$, then the generalized coarse
index $\textindsmall D$ vanishes in $K_0(C^\ast_G(M))$.

\medskip\noindent
{\bf Proposition 2:} Let
$\widetilde D$ the lift of a generalized elliptic operator in
$L^2(\widetilde M,S)$. 
In the ungraded case, 
if there is a gap in the spectrum of $\widetilde D$, then the index
$\textindsmall D$ vanishes in $K_1(\cgm)$.

\medskip\noindent
{\bf Corollary:} Let $M$ be a complete spin manifold. If
$M$ has a metric of uniformly positive scalar curvature in some
coarse class, then
the generalized coarse index of the spinor Dirac operator vanishes.

We now embark on the task of computing the $K$-theory of this 
algebra and of coarse indices.

Let $(M,d)$ be a proper metric space. For any subset $U\subset M$ and
$R>0$, we denote by $\textpen(U,R)$ the open neighborhood of $U$
consisting of points $x\in M$ for which $d(x,U)<R$. Let $A$
and $B$ be closed subspaces of $M$ with $M=A\cup B$. We then
say that the decomposition  $(A,B)$ is a {\it coarsely excisive pair}
if for each $R>0$ there is an $S>0$ such that $$\textpen(A,R)
\cap \textpen(B,R)\subseteq \textpen(A\cap B,S).$$
We wish to analyze this decomposition in the following context.

Given general $C^\ast$-algebras $\cala$, $\calb$ and $\calm$ for which
$\calm=\cala+\calb$, we have the Mayer-Vietoris sequence
$$\cdots\longrightarrow K_{j+1}(\,\calm\,)\longrightarrow
K_j(\,\cala\cap\calb\,)\longrightarrow K_j(\,\cala\,)\oplus
K_j(\,\calb\,)\longrightarrow
K_j(\,\calm\,) \longrightarrow\cdots$$
The standard proof for the existence of such a sequence is
developed from the isomorphism $K_\ast({\mathcal T})\cong K_{\ast-1}(\calm)$,
where ${\mathcal T}$ is the suspension of $\calm$. A short discussion of this
construction is given in \cite{HigsonRoeYu}.
We are in particular interested in exploiting the boundary map
$\partial: K_j(\,\calm\,)\to K_{j-1}(\,\cala\cap\calb\,)$ to
transfer information about the index of the Dirac operator
on a complete noncompact manifold $M$ to information about that
on some hypersurface $V$. For our purposes, we wish to set $\calm$
to be the generalized Roe algebra $\cgm$ on $M$, while $\cala$
and $\calb$ represent analogous operator algebras on closed subsets
$A$ and $B$, where $(A,B)$ form a coarsely excisive decomposition
of $M$. To construct the boundary map in question, we require
a few technical lemmas and notion of equivariant operators with
generalized finite propagation on a subset of $M$. The proof of the
first lemma follows the same argument as that in \cite{HigsonRoeYu}
and is stated without proof.

\medskip\noindent
{\bf Definition:} Let $A$ be a closed subspace of a proper metric
space $M$. Denote by $D^\ast_G(A,M)$ the $C^\ast$-algebra
of all operators $T$ in $\dgm$ 
such that $\textsupp T\subseteq
\textpen(\pi^{-1}(A), R)\times \textpen(\pi^{-1}(A), R)$,
for some $R>0$. Let $C^\ast_G(A,M)$ be the quotient $D^\ast_G(A,M)/I$.

\medskip\noindent
{\bf Lemma 2:} Let $(A,B)$ be a decomposition of
$M$. Then
\begin{enumerate}
\item $C^\ast_G(A,M)+C^\ast_G(B,M)=\cgm$.
\item $C^\ast_G(A,M)\cap C^\ast_G(B,M)=C^\ast_G(A\cap B, M)$ if
in addition we assume that $(A,B)$ is coarsely excisive.
\end{enumerate}

\medskip\noindent
{\bf Lemma 3:} Suppose $V\subset M$ and $\pi_1(V)$ injects into
$\pi_1(M)$. There is an isomorphism $K_\ast(C^\ast_G(V))\cong
K_\ast(C^\ast_G(V,M))$.

\smallskip\noindent
{\it Proof:} Let $\pi:\widetilde M\to M$ be the projection map.
Consider the $C^\ast$-algebra  $C^\ast(\textpen(\pi^{-1}(V),n),\pi_1(M))$
given by the quotient by $I$ of the $C^\ast$-algebra of
locally compact, $\pi_1(M)$-equivariant operators on
the $n$-neighborhood penumbra $\textpen(\pi^{-1}(V),n)$.
Then $$C^\ast_G(V,M) =\lim_{\longrightarrow}
C^\ast(\textpen(\pi^{-1}(V),n),\pi_1(M)).$$
The inclusion map $i:\pi^{-1}(V)\hookrightarrow
\textpen(\pi^{-1}(V),n)$ is a coarse equivalence.
Since by the construction the generalized Roe algebra its
operators are defined up to their bounded parts, the map $i$
induces a series of isomorphisms
$$
\begin{array}{ccl}
K_\ast(C^\ast(\pi^{-1}(V),\pi_1(M)))&\cong & 
K_\ast(C^\ast(\textpen(\pi^{-1}(V),n),\pi_1(M))\\
&\cong& 
  K_\ast(\,\,\lim
     C^\ast(\textpen(\pi^{-1}(V),n),\pi_1(M))\\
 &\cong&K_\ast(C^\ast_G(V,M)).
\end{array}
$$
Since $\pi_1(V)\hookrightarrow\pi_1(M)$ is an injection, the inverse
image $\pi^{-1}(V)\subseteq \widetilde M$ is a disjoint union
of isomorphic copies of $\widetilde V$, parametrized by the
coset space $\pi_1(M)/\pi_1(V)$. Therefore, there
is a one-to-one correspondence between $\pi_1(M)$-equivariant
operators on $\pi^{-1}(V)$ and $\pi_1(V)$-equivariant operators
on $\widetilde V$. Hence $C^\ast(\pi^{-1}(V),\pi_1(M))\cong
C^\ast_G(V)$. We then have $K_\ast(C^\ast_G(V))\cong
K_\ast(C^\ast_G(V,M))$, as desired. \hfill $\Box$

\medskip
Let $(A,B)$ be a coarsely excisive decomposition of $M$ such
that $V=A\cap B$ satisfies $\pi_1(V)\hookrightarrow\pi_1(M)$.
The boundary operator 
$\partial: K_j(C^\ast_G(A,M)+C^\ast_G(B,M))\longrightarrow
K_{j-1}(C^\ast_G(A,M)\cap C^\ast_G(B,M))$
arising from the coarse Mayer-Vietoris sequence is by the previous
lemmas truly a 
map $$\partial: K_\ast(C^\ast_G(M))\to K_{\ast-1}(C^\ast_G(V)).$$

\medskip\noindent
{\bf Theorem: (Boundary of Dirac is Dirac)}
Consider a coarsely excisive decomposition $(A,B)$ of $M$ and
let $V=A\cap B$. If $\partial : K_\ast(\cgm)\to K_{\ast-1}(C^\ast_G(V))$
is the boundary map from the Mayer-Vietoris sequence derived above,
then we have $\partial\,(\textind_M(D))=\textind_V(D)$.

\medskip\noindent
{\bf Remark:} Here $\textind_M(D)$ and $\textind_V(D)$ represent
the generalized coarse indices of the spinor Dirac operators on
$M$ and $V$, respectively. We will continue to use a subscript
if the space to which the index is related is ambiguous. 
The ``boundary of Dirac is
Dirac" principle is essentially equivalent to Bott periodicity
in topological $K$-theory. In all cases considered here, there
are commutative diagrams relating topological boundary to the
boundary operator arising in the $K$-theory of $C^\ast$-algebras, and
on the topological side, a consideration of symbols suffices.
See \cite{RoePartitioning}, 
\cite{HigsonDuality}, \cite{RoeIndex} and \cite{Wu}. 

\medskip\noindent
{\bf Theorem 1:} The $n$-fold product $M$ of punctured two-dimensional
tori  does not
have a metric of uniform positive scalar curvature in the same coarse
equivalence class as the positive hyperoctant with its
standard Euclidean metric.

\smallskip\noindent
{\it Proof:}
Consider the projection map
$p: M\to \reals_{\ge0}^n$ from the multifold product
$M=\,\ptorusn$ to the positive hyperoctant, where each
component $p_i$ is the quasi-isometric projection of
the punctured torus onto the positive reals numbers.
Take a hypersurface $S\subset  \reals_{\ge0}^n$ sufficiently
far from the origin so that the inverse image of every point
on $S$ is an $n$-torus, and so the space $V$ is coarsely equivalent 
to the $(2n-1)$-dimensional noncompact manifold $\reals^{n-1}\times \,T^n$. 
The complement of the hypersurface $V$ consists of two noncompact components.
Define $A$ to be the closure of the
component containing the inverse image $p^{-1}
(\mathbf{0})$ of the origin in  $\reals^n_{\ge0}$. 
Take $B$ the closure of $M\backslash A^\circ$. 
Then the pair $(A,B)$ forms a coarsely excisive decomposition 
of the space $M$ whose intersection is $A\cap B=V$.

\vskip -.4in
\hskip .5in
\begin{figure}[ht]
\begin{tabbing}
%\vskip -1cm
%\hskip 5cm
\hskip 3cm\=
\setlength{\unitlength}{.75cm}
%\begin{picture}(4,6)
%\put(0,1){\vector(1,0){3.5}}
%\put(0,1){\vector(-1,0){3.5}}
%\put(0,0){\vector(0,1){5}}
%\qbezier(-2,4.5)(0,.75)(2,4.5)
%\put(2.3,4){\makebox(0,0)[bl]{{\footnotesize $y=\sqrt{9x^2+1}$}}}
%\end{picture}~~
%\=\hskip 3cm
%\vskip -1cm
\hskip 1cm
\setlength{\unitlength}{.75cm}
\begin{picture}(4,6)
\put(0,0){\vector(1,0){5}}
\put(0,0){\vector(0,1){4.5}}
\qbezier(2.5,4.5)(0,0)(4.5,2.5)
\put(4.5,1.8){\makebox(0,0)[bl]{{\footnotesize $S$}}}
\end{picture}
\end{tabbing}
\caption{The hypersurface $S$ in $\reals^n_{\ge0}$.}
\end{figure}

\smallskip
Consider the generalized coarse index $\textind_M(D) \in K_\ast(\cgm)$
of the lifted classical Dirac operator on the pullback spinor bundle
of the universal cover $\widetilde M$. Note that $\pi_1(M)$
is the $n$-fold product $F_2\times\cdots\times F_2$ of free groups, and
that $\pi_1(V)\cong\pi_1(\reals^{n-1}\times T^n)
\cong\integers^n$. Hence there is an injection
$\pi_1(V)\hookrightarrow\pi_1(M)$ and the $K$-theoretic Mayer-Vietoris
sequence applies. The boundary map $\partial$ of
this sequence satisfies 
$\partial\,(\textind_M(D))=\textind_V(D)\in K_\ast(C^\ast_G(V))$. 
However, $V$ is coarsely equivalent to the hypersurface
$\reals^{n-1}\times T^n$, so the index $\textind_V(D)$ can be taken
to live in $K_{\ast-1}(C^\ast_G(\reals^{n-1}\times T^n))$. Note that
$n$ will be taken to be at least $2$. There is yet
another boundary map $K_{\ast-1}(C^\ast_G(\reals^{n-1}\times T^n))
\to K_{\ast-n+1}(C^\ast_G(\reals\times T^n))$ by peeling off $n-2$
copies of the real line. This boundary map (or composition of
$n-2$ boundary maps) preserves index. 

Recall that  $D_G^\ast(M)$ is the norm closure of the $C^\ast$-algebra
of all locally compact, $\pi_1(M)$-equivariant, generalized finite
propagation operators on $L^2(\widetilde M)$, and $I\subseteq
D_G^\ast(M)$ is the closure of the
ideal of such operators $T$ that satisfies the condition that
$(\pi\times\pi)(\textsupp T)$ is bounded in $M\times M$.
The short exact sequence $0\longrightarrow I\longrightarrow 
D^\ast_G(\reals\times T^n)
\longrightarrow C^\ast_G(\reals\times T^n)
\longrightarrow 0$ gives rise to the
six-term exact sequence in $K$-theory:
$$\diagram
  K_0(I)\rto& K_0(D^\ast_G(\reals\times T^n))\rto& 
	     K_0(C^\ast_G(\reals\times T^n))\dto\\
	K_1(C^\ast_G(\reals\times T^n))
	 \uto& K_1(D^\ast_G(\reals\times T^n))\lto& 
		    K_1(I)\lto
		    \enddiagram $$
Notice that the map $K_\ast(I)\to K_\ast(D^\ast_G(\reals\times T^n))$ 
induced by the inclusion is the zero map by an Eilenberg swindle
argument. 
Hence both maps $K_\ast(D^\ast_G(\reals\times T^n))\to
K_\ast(C^\ast_G(\reals\times T^n))$ are injections. 

If $n$ is even, 
the generalized coarse index $\textind_{\reals\times T^n}(D)$ of 
the Dirac operator $D$ resides in $K_1(D^\ast_G(\reals \times T^n))$. 
 Certainly the image of this  index under
the boundary map $K_1(D^\ast_G(\reals\times T^n))\to
K_0(D^\ast_G(T^n))$ is the index $\textind_{T^n}(D)$ of $D$ on the
$n$-torus. Since $T^n$ does not have a metric of positive scalar
curvature at all, the obstruction $\alpha(T^n,f)\in 
K_\ast(C^\ast(\integers^n))$, where $f: T^n\to B\integers^n$ is 
the classifying map, is nonvanishing. This ``index" is constructed
by Rosenberg in \cite{Rosenberg2}. This special case of the
Gromov-Lawson-Rosenberg conjecture holds for the group $\integers^n$.
This index maps to our generalized coarse index $\textindsmall_{T^n} D$
under the isomorphism
$K_\ast(C^\ast(\integers^n))
\cong K_\ast(D^\ast_G(T^n))$.
Hence the index of $D$ in  $K_1(D^\ast_G(\reals
\times T^n))$ is nonzero, and its projection onto the group
$K_1(C^\ast_G(\reals\times T^n))$ is nonzero as well. This argument
gives us the necessary index obstruction.  

If $n$ is odd, we apply the same argument as above with respect to the
map $K_0(D^\ast_G(\reals\times T^n))\to
K_0(C^\ast_G(\reals\times T^n))$. \hfill $\Box$

\bigskip
The extension of this method to multifold products of 
hyperbolic manifolds involves the Margulis lemma, which
states that in such a space there exists a small positive constant
$\mu=\mu_n$ such that the subgroup $\Gamma_\mu(V, v)\subseteq
\pi_1(V,v)$ generated by loops of length less than or equal to
$\mu$ based at $v\in V$ is almost nilpotent, i.e. it contains
a nilpotent subgroup of finite index.
It can be shown that there exist such
cusps, or submanifolds $C\subset V$ with compact convex boundary
containing $v$, such that $C$ is diffeomorphic to the product
$\partial C\times\reals_+$, where $\partial C$ is diffeomorphic
to an $(n-1)$-dimensional nilmanifold with fundamental
group containing $\Gamma_\mu(V,v)$. Here a nilmanifold signifies
a quotient $N/\Gamma$ of a nilpotent Lie group by a cocompact
lattice $\Gamma$. The nilmanifolds that arise in this context
as boundaries of pseudospheres will have a naturally flat
structure.

\medskip\noindent
{\bf Theorem 2:} An $n$-fold product of hyperbolic
manifolds has no uniform positive scalar curvature metric coarsely equivalent 
to the usual Euclidean metric on the
positive Euclidean hyperoctant.

\smallskip\noindent
{\it Proof:} Without loss of generality, it suffices to
consider the case when the noncompact hyperbolic spaces
have only one cusp.  Let $m$ be the dimension of this
product manifold. As in the multifold product of tori,
there is a positive $b\in\reals$ such that 
on each hyperbolic space $\calh_i$ the
inverse image of each point $x\ge b$ under
the projection $\calh_i\to\reals_{\ge0}$ is by Margulis'
lemma a flat compact
connected Riemannian manifold of finite dimension. Consider
the inverse image $\calv$ under the induced product map $p:\calh_1
\times\cdots\times\calh_m
\to\reals^m_{\ge0}$ of the same hypersurface as described in
the previous theorem. By Bieberbach's theorem, every flat
compact connected Riemannian manifold admits a normal Riemannian
covering by a flat torus of the same dimension. Hence $\calv$ is covered
by some product of Euclidean space and a higher-dimensional torus.
Any metric of positive scalar curvature on $\calv$ would certainly
lift to such a metric in this covering space. Using the same induction argument as before, we show that such a metric is obstructed by the
presence of a nonzero Dirac class. \hfill\ $\Box$

\bigskip\bigskip\noindent
III. Noncompact Quotients of Symmetric Spaces: A Special Case

%\bigskip\noindent
%{\bf Lemma:} \cite{Wolf} The double coset space $\slnz\backslash\slnr/\sonr$
%parametrizes the isometry classes of flat Riemannian $n$-tori
%KK unit volume.

\bigskip
The Iwasawa decomposition gives a unique way of expressing
the group $\slnr$ as a product $\slnr = NAK$, where
$N$ is the subgroup of standard unipotent matrices (upper
triangular matrices with all diagonal entries equal to
$1$), $A$ the subgroup of $\slnr$ consisting of 
diagonal matrices with positive entries, and $K$ the
orthogonal subgroup $\sonr$. 
The effect of taking the double quotient of $\slnr$ by
both $\sonr$ and $\slnz$ on the Iwasawa decomposition
is that we are left with classes of matrices 
represented by those of the form $n^\ast a^\ast$, 
where $n^\ast$ is represented by a unipotent matrix
and $a^\ast= \textdiag(a_1,\ldots, a_n)$ is
diagonal with weakly increasing entries
$a_1\le a_2\le \cdots\le a_n$. Let $N^\ast$ be the
iterated circle bundle that arises upon taking a quotient
of $N$ by $\slnz$ (for example, when $n=3$, the group
$N$ is a Heisenberg group). Let $A^\ast$ be the
semigroup of matrices (note that it is not closed under
inversion) with such an increasing condition on
the entries.  Consider the map $\rho:
\slnz\backslash \slnr/\sonr \to A^\ast$ 
given by $n^\ast a^\ast\mapsto a^\ast$. This map denotes a fiber bundle
with fiber $N^\ast$ over every point $a^\ast\in A^\ast$. Notice
that the Iwasawa decomposition gives another way of observing
the dimension of $\slnr/\sonr$, the sum of $n-1$ dimensions from $A$ and
$\frac{n(n-1)}{2}$ from $N$. The quotient $\slnz\backslash \slnr/\sonr$
then has $n-1$ noncompact directions coming from $A^\ast$ and 
$\frac{n(n-1)}{2}$ compact directions from $N^\ast$.

The space $A^\ast$ is identifiable with the subset of $(n-1)$-dimensional
Euclidean space given by
$\{(a_1,\ldots, a_n): 0<a_1\le\cdots\le a_n, ~a_1\cdots a_n=1\}$.
We wish to construct a hypersurface in $A^\ast$ whose inverse
image under $\rho$ is an iterated circle bundle over Euclidean
space. We recall that a geodesic in the quotient space  $\slnr/\sonr$ is of the
form $t\mapsto e^{t\Lambda}\cdot\sonr$, where $\Lambda$ is an $n\times n$
symmetric matrix with zero trace. We shall construct an appropriate
hypersurface in $A^\ast$ by taking the union of sufficiently many geodesics.

Consider a geodesic $t\mapsto e^{t\Lambda}\cdot\sonr$ in the
quotient $\slnr/\sonr$ where
$\Lambda=\textdiag(\lambda_1,\ldots, \lambda_n)$ with $\lambda_1\le
\cdots\le\lambda_n$ and $\lambda_1+\cdots
+\lambda_n=0$. Then the geodesic is a map $t\mapsto \textdiag(e^{\lambda_1t},
\ldots, e^{\lambda_nt})\cdot\sonr$. Because two symmetric traceless
matrices of the form $\Lambda_1=\textdiag(\lambda_1,\ldots, \lambda_n)$
and $\Lambda_2=\textdiag(\alpha\lambda_1,\ldots, \alpha\lambda_n)$ give
the same geodesic image for $\alpha>0$,
we can normalize the $\lambda$-vector so
that $\lambda_1=-1$. In addition, let $\mu_2=1$ and
$\mu_i\in [\,i-1,i\,]$ for
$i=3,\ldots, n$ and let $m=\sum_{i=2}^n\mu_i$. Set $\lambda_i = \mu_i/m$.
Then $\Lambda=\textdiag(\lambda_1,\ldots, \lambda_n)$ is symmetric 
and traceless with weakly increasing entries.

\medskip\noindent
{\bf Lemma 4:} Each ordered $(n-1)$-tuple $(\mu_2,\ldots, \mu_n)
\in [\,2,3\,]\times\cdots\times [\,n,n+1\,]$ gives rise to a unique
geodesic $t\mapsto e^{t\Lambda}\cdot\sonr$, up to reparametrization .

\smallskip\noindent
{\it Proof:} Suppose  $(\mu_2^{(1)},\ldots, \mu_n^{(1)}), 
(\mu_2^{(2)},\ldots, \mu_n^{(2)})
\in [\,2,3\,]\times\cdots\times [\,n,n+1\,]$ give rise to the same
geodesic. The two vectors correspond to the traceless 
matrices $\Lambda_1=\textdiag(-1,\lambda_2^{(1)},\ldots, \lambda_n^{(1)})$
and $\Lambda_2=\textdiag(-1,\lambda_2^{(2)},\ldots, \lambda_n^{(2)})$, 
respectively.  Let $\nu_1=(-1,\lambda_2^{(1)},\ldots, \lambda_n^{(1)})$
and $\nu_2=(-1,\lambda_2^{(2)},\ldots, \lambda_n^{(2)})$. Obtain
the normalized matrices $\Lambda_1^\ast$ and $\Lambda_2^\ast$ by
dividing the entries by the Euclidean norms $||\nu_1||$ and
$||\nu_2||$. By assumption, $\Lambda_1^\ast =\Lambda_2^\ast$;
in other words, $\frac{\nu_1}{||\nu_1||}= \frac{\nu_2}{||\nu_2||}$.
But then $\theta = \cos^{-1}\frac{\nu_1\cdot\nu_2}{||\nu_1||||\nu_2||}
=0$, so $\nu_1$ and $\nu_2$ are parallel. Since their first coordinates
coincide, they are identical. Hence
$\frac{(\mu_2^{(1)},\ldots,\mu_n^{(1)})}{\mu_2^{(1)}+\cdots+\mu_n^{(1)}}=
\frac{(\mu_2^{(2)},\ldots,\mu_n^{(2)})}{\mu_2^{(2)}+\cdots+\mu_n^{(2)}}$.
Since each $\mu_i^{(j)}$ is positive, it can be written as the
square of some other positive number. Arguing as before, we see that
$(\mu_2^{(1)},\ldots, \mu_n^{(1)})= \beta(\mu_2^{(2)},\ldots, \mu_n^{(2)})$. 
for some $\beta>0$. Since $\mu_2^{(1)}$ and $\mu_2^{(2)}$ both
equal one, the vectors are coincident. \hfill $\Box$

\medskip
The space $A^\ast$ of diagonal matrices in $\slnr$ with weakly
increasing entries is itself simply-connected of dimension $n-1$,
and its space at infinity is an $(n-2)$-dimensional simplex $P$. Let $W$
be the union of geodesics constructed above, with
$(\mu_2,\ldots, \mu_n)$ ranging in the product
$[\,2,3\,]\times\cdots\times [\,n,n+1\,]$ of intervals. The subset
$W' \subset A^\ast$ is an $(n-1)$-dimensional space with boundary
$V'=\partial W'$. The $(n-2)$-dimensional
hypersurface $V'$ is also simply-connected
whose space at infinity is homeomorphic to an $(n-3)$-sphere
that is disjoint with $\partial P$. The space $V'$ itself
is coarsely equivalent to $\reals^{n-2}$.

\medskip\noindent
{\bf Theorem 3:} The double quotient group $\slnz\backslash \slnr/\sonr$
does not have a uniform positive scalar curvature metric that is coarsely
equivalent to the natural one inherited from $\slnr$.

\smallskip\noindent
{\it Proof:} Let $M= \slnz\backslash \slnr/\sonr$ and recall
the projection map $\rho: M\to A^\ast$
given by $n^\ast a^\ast\mapsto a^\ast$. Since the fiber over each
point is an arithmetic quotient of the group of unipotent matrices,
the inverse image $V\equiv\rho^{-1}(V')$ is coarsely equivalent to an iterated
circle bundle over Euclidean space. Moreover, $V$ partitions the
space into a coarsely excisive pair whose closures $(A,B)$ satisfy
the equalities $A\cup B= M$ and $A\cap B=V'$ (note that $B$ can be taken
as $\rho^{-1}(W')$ and $A$ the closure of its complement). 
If $\textind_M(D)$ denotes the
generalized coarse index of the classical
spinor Dirac operator on $\widetilde M$, then
the Mayer-Vietoris map $\partial: K_\ast(C^\ast_G(M))
\to K_{\ast-1}(C^\ast_G(V))$ defined in the previous chapter 
satisfies $\partial\,(\textind_M(D))=\textind_{V'}(D)=
\textind_{\reals^{n-2}\times U^m}(D)$, where $U^m$ is the compact fiber
of the iterated circle bundle of dimension $m=\frac{n(n-1)}{2}$.
Applying the same argument as before, it suffices to show
that the index of the Dirac operator in $K_\ast(D^\ast_G(U^m))$
is nonzero. However, $U^m$ is a quotient of a nilpotent 
group with a cocompact lattice, and hence by Gromov
and Lawson \cite{GromovLawson2} has no
metric of positive scalar curvature at all. As with the theorem
for punctured tori, there is a nonvanishing Rosenberg index 
$\alpha(U^m)\in K_\ast(C^\ast(\pi))$, where $\pi=\pi_1(U^m)$,
which maps to the generalized coarse index in $K_\ast(D^\ast_G(U^m))$,
as desired. Here the Gromov-Lawson-Rosenberg conjecture is
true since $U^m$ is a nilmanifold. \hfill $\Box$

\bigskip\bigskip\noindent
IV. The General Noncompact Arithmetic Case

\bigskip
To understand how we might achieve a similar result for general
double quotient spaces $\Gamma\backslash G/K$, we appeal
to the following. 

\bigskip\noindent
{\bf Picture from Reduction Theory:} There is a compact polyhedron $P$ and
a Lipschitz map $\pi : M\to  cP$, where $cP$ is the open cone on  $P$
so that (1) every point inverse deform retracts to an
arithmetic manifold, (2)  $\pi$
respects the radial direction, and (3) all point inverses have
uniformly bounded size.

\bigskip
        Indeed, the polyhedron $P$ is the geometric realization 
of the category of proper $\rationals$-parabolic
subgroups of $G$, modulo the action of $\Gamma$, and  $\pi^{-1}$ of the
barycenter of a simplex is the arithmetic symmetric space
associated to that parabolic.  Concretely, for $\hbox{SL}_n(\integers)
\subset \hbox{SL}_n(\reals)$, the space $P$ is an $n-2$
simplex, the parabolics correspond to flags, and the associated
arithmetic groups have a unipotent normal subgroup with quotient equal
to a product of $\hbox{SL}_{m_i}(\integers)$, where the $m_i$
 are sizes of the blocks occurring in
the flag.  As one goes to infinity, the unipotent directions shrink in
diameter and are responsible for the finite volume property of the
lattice quotient, while the other parabolic directions remain of
bounded size.

Alternatively, for any choice of basepoint in the
homogeneous space, there are constants $C$ and $D$ that satisfy the
following condition: if $x$ is a given point and
$Q_x$ is the largest parabolic subgroup associated with a simplex
whose cone contains $x$ within its $C$-neighborhood, then
the orbit of $x$ under  $Q_x$ has diameter less than $D$.
Note the empty simplex means that there is a compact
core which is stabilized by the whole group.

        This picture can essentially be ascertained from \cite{BorelSerre},
\cite{Saper}; the
fact that $K\backslash G/\Gamma$ has finite Gromov-Hausdorff distance from
$cP$ is asserted in \cite{Gromov}. However,
one needs a key estimate about the ``coarse isotropy." Details are given
in unpublished work of Eskin \cite{Eskin}; some are given below.

As a guide the reader should think through the picture
suggested by a product of  hyperbolic manifolds.  Each
hyperbolic manifold contributes to $cP$
either a point, in the compact case,  or the
open cone on a finite set of points, in the case of cusps.  Thus $P$
is a join of some number of finite sets.  Using this model, we find
that the
inverse image of any point in the interior of any simplex is exactly a
product of closed hyperbolic manifolds, cores of hyperbolic
manifolds, and flat manifolds.

\medskip
Let us now consider the unique decomposition \cite{Eberlein}
of a semisimple Lie group $G=N_x\cdot
A_x\cdot K$, where $x$ is a point on the space of infinity of $G/K$.
If $\Gamma$ is an arithmetic lattice in $G$, we
are interested in knowing how $\Gamma$ acts on $G/K$. In other words,
we ask how $\Gamma$ acts on this particular coordinate system. Let
$g=nak\in G$. If $\gamma\in\Gamma$, let $\gamma g=n'a'k'$. Notice
that $N_x\cap \Gamma$ is a lattice in $N_x$, and acts cocompactly
on $N_x$. Consider the projection $\rho: \Gamma\backslash G/K\to A_x^\Gamma$
given by $n^\ast a^\ast \mapsto a^\ast$, where
$A_x^\Gamma$ is a fundamental domain of $A_x$ under
the action of $\Gamma$. The fiber above each point is
a compact manifold arising from the action of $N_x\cap \Gamma$ on $N_x$.
Let $\frakf+\frakp$ be the usual Cartan decomposition and let
$X\in\frakp$ be the element such that $dp(X)=\gamma_{px}'(0)$.
If $Z(X)$ is defined by $Z(X)=\{Y\in\frakg: [\,X,Y\,]=0\}$, then
$\fraka=Z(X)\cap \frakp$ is the unique maximal abelian subspace of
$\frakp$ that contains $X$. By definition $A_x=\exp\,(\fraka)$.

We are interested in constructing a hypersurface $V'$ in the 
space $A_x^\Gamma$
whose inverse image $V\equiv\rho^{-1}(V')$ under the projection map $\rho$
provides us with an appropriate excisive decomposition of $\Gamma
\backslash G/K$ for which the $K$-theoretic Mayer-Vietoris sequence
is applicable. Consider the chamber decomposition of $\fraka$.
%\hskip 4in
%\begin{figure}[h]
%%(means put here or at bottom)
%%\begin{figure}[b]
%\begin{center}
%  \includegraphics[height=2.3in,width=8in,angle=-90,keepaspectratio]{eskin.eps}
%  \end{center}
%\caption{The Weyl chambers of an maximal abelian subspace $\fraka\subseteq
%\frakg$.}
%\end{figure}
The Weyl group $W$ acts on these chambers via the hyperplanes. Consider
the Bruhat decomposition $G= \coprod_{w\in W}BwB$ of $G$, and let
$\gamma\in BwB$ for some $w\in W$. Recall that, if $g=nak$, 
we write $\gamma g=n'a'k'$. Denote by $\sum^+$ the set of
positive roots of $\fraka^\ast$ and $\sum^-$ the set of negative roots. Let
$\calr= \sum^-\cap \,w\sum^+$ be the set of roots that start
out positive but are made negative under the action of $w$.
For some positive
reals constants $c_\alpha$, the following equation holds 
\cite{Abels},\cite{Eskin}:
$$a' = wa-\sum_{\alpha\in\calr}c_\alpha\alpha(a) +~ O(1),$$
where $a$ and $a'$ are viewed as elements of the Lie algebra $\fraka$.
The implications of this equation are as follows. Consider
an element $a$ in the positive Weyl chamber $\calc(a)$. 
Then the intersection $\Gamma(a)\,\cap\, \calc(a)$ of the orbit $\Gamma(a)$ of
$a$ under $\Gamma$ and the Weyl chamber $\calc(a)$ containing $a$
has a bounded diameter, uniformly in $a$. In other words, if $\gamma(a)$ stays
in the same positive Weyl chamber, then $w$ is the identity and
$\calr$ is empty. Hence $a'=a+o(1)$, implying that $a'$ 
can be found at a uniformly
bounded distance from $a$ itself. In this event, the action of $\gamma$
corresponds to a translation of $a$
to (possibly) the compact fiber direction of $\Gamma\backslash G/K$.
It is also a general fact that the action of any $g\in\Gamma$ 
will take the vertex of any subsector (as drawn in
the previous figure) to the vertex of an analogous subsector.
With this machinery, we are able to prove the following.

\medskip\noindent
{\bf Theorem 4:} The double quotient space $M=\Gamma\backslash G/K$
has no metric of uniform positive scalar curvature in the same coarse
class as the natural metric inherited from $G$.

\smallskip\noindent
{\it Proof:} 
The picture from reduction theory provides
a polyhedron $P$ in the space of infinity of the positive Weyl chamber and
an open cone $W=cP$ on $P$ in
$\calc^+$ oriented so that the distance from $W$ to any hyperplane 
$\alpha=0$ will exceed the quantity
$\sup_{a\in \calc^+}\textdiam (\Gamma(a)\cap \calc^+)$. This quantity
is finite by \cite{Eskin}, \cite{Abels}.
Let $V'=\partial W$ and $V=\rho^{-1}(V')$
in $M=\Gamma\backslash G/K$. Then $V$ is a hypersurface that
induces a decomposition $(A,B)$
of $M$. The fundamental group $\pi_1(V)=N\cap \Gamma$
injects into $\pi_1(\Gamma\backslash G/K)=\Gamma$, satisfying
the requirement needed in the construction of the Mayer-Vietoris
sequence. In the most general case, the space $V$ is coarsely equivalent 
to a bundle over Euclidean space whose fiber consists of two components:
a nilmanifold $N$ and (possibly) a compact homogeneous manifold $H$. 
In the absence of such an $H$, the argument follows exactly as it does for
$\slnz\backslash\slnr/\sonr$. In the presence of a compact homogeneous
manifold, we may pass the coarse index of the Dirac operator to $\reals\times H$ 
and use the usual Rosenberg obstruction on $H$ as in Theorem 1 to
obtain our desired result. \hfill $\Box$

\end{document}